\documentclass{amsart}

\usepackage{amssymb}
\usepackage{xcolor}
\usepackage{amsmath}
\usepackage{bm}

\newcommand{\half}{\frac{1}{2}}
\newcommand{\thalf}{\tfrac{1}{2}}
\newcommand{\summ}{\mathop{{\sum}^{\star}}}

\newcommand{\sym}{{\rm sym}}

\makeatletter
\@namedef{subjclassname@2020}{\textup{2020} Mathematics Subject Classification}
\makeatother
\numberwithin{equation}{section}

\newtheorem{theorem}{Theorem}[section]

\newtheorem{lemma}[theorem]{Lemma}

\begin{document}

\title{On the error term in a mixed moment of $L$-functions
}

\author{Rizwanur Khan and Zeyuan Zhang}

%\address{
%Department of Mathematical Sciences\\ The University of Texas at Dallas\\ 800 W. Campbell Rd.\\ Richardson, TX 75080-3021}
%\email{rk2357@gmail.com}

\address{
Department of Mathematics\\ University of Mississippi\\ University, MS 38677}
\email{rrkhan@olemiss.edu, zzhang8@go.olemiss.edu}

\thanks{The authors were supported by NSF grant DMS-2001183. The first author was also supported by NSF grant DMS-2140604 and the Simons Foundation (award 630985). }

\begin{abstract} There has recently been some interest in optimizing the error term in the asymptotic for the fourth moment of Dirichlet $L$-functions and a closely related mixed moment of $L$-functions involving automorphic $L$-functions twisted by Dirichlet characters. We obtain an improvement for the error term of the latter.
\end{abstract}

\maketitle

\section{Introduction}

Heath-Brown \cite{hb} was the first to prove an asymptotic formula for the fourth moment of the Riemann Zeta function on the critical line with a power saving error term. The analogous asymptotic for Dirichlet $L$-functions required new ideas and came much later. Young \cite{you} proved that
\begin{align}
\label{4thasymp} \frac{1}{\varphi^*(q)} \summ_{\chi\bmod q} |L(\thalf, \chi)|^4=P_4(\log q)+O(q^{-\frac{5}{512}+\epsilon }),
\end{align}
where $q$ is prime, $\summ$ restricts the sum to primitive Dirichlet characters, $\varphi^*(q)$ is the number of primitive Dirichlet characters, and $P_4$ is a degree four polynomial. Wu \cite{wu} obtained a power saving asymptotic for general moduli $q$ with an error term of $O(q^{-\frac{11}{448}})$.  Blomer, Fouvry, Kowalski, Michel, and Mili\'cevi\'c \cite{bfkmm1} revisited the problem, and introducing new ideas, they significantly sharpened Young's error term for $q$ prime to $O(q^{-\frac{1}{32}+\epsilon})$. %,and to $O(q^{-\frac{1}{24}+\epsilon})$ assuming the Ramanujan conjecture.
Later Blomer et al. \cite{bfkmm2} returned to this problem once more and used an idea of Shparlinski and T. Zhang \cite{shpzha} to reduce the error term even further to $O(q^{-\frac{1}{20}+\epsilon})$, and to $O(q^{-\frac{1}{14}+\epsilon})$ assuming the Ramanujan Conjecture.\footnote{Actually the weaker bound $O(q^{-\frac{1}{16}+\epsilon})$ is claimed on the Ramanujan conjecture at the end of section 4 in \cite{bfkmm2}, but this seems to be a calculation mistake. This has also been observed in \cite[Theorem 1.1]{zac} with $\theta=0$ and $L=1$.}

In \cite{bfkmm1}, Blomer et al. also considered a related problem which is one step greater in difficulty. They established the mixed moment asymptotic
\begin{align}
\label{mixed-original} \frac{1}{\varphi^*(q)} \summ_{\chi\bmod q} L(\thalf, f\otimes \chi)\overline{L(\thalf, \chi)}^2=\frac{L(1,f)^2}{\zeta(2)}+O(q^{-\frac{1}{68}+\epsilon}),
\end{align}
for prime $q$, where $f$ is a fixed (holomorphic or non-holomorphic) Hecke eigenform for $SL_2(\mathbb{Z})$. There has been some interest in improving the error term in this asymptotic as well. Shparlinski \cite{shp} improved the error term in \eqref{mixed-original} slightly to $O(q^{-\frac{1}{64}+\epsilon})$. The goal of this paper is to make a more significant improvement to this error term. In fact when $f$ is holomorphic, the quality of our error term is almost as good as the error term in the fourth moment of Dirichlet $L$-functions. This is achieved by an idea which simplifies previous approaches to the problem in a critical range. 
\begin{theorem}\label{main} Let $q$ be an odd prime. Let $\eta_f=\frac{1}{22}$ for $f$ a holomorphic Hecke-cusp form, and $\eta_f=\frac{5}{152}$ for $f$ a Hecke Maass-cusp form, for $SL_2(\mathbb{Z})$. We have
 \begin{align}
\label{main-asymp} \frac{1}{\varphi^*(q)} \summ_{\chi\bmod q} L(\thalf, f\otimes \chi)\overline{L(\thalf, \chi)}^2=\frac{L(1,f)^2}{\zeta(2)}+O(q^{-\eta_f+\epsilon}).
\end{align}
\end{theorem}

We also mention that the more difficult related asymptotic
\begin{align}
\label{hard} \frac{1}{\varphi^*(q)} \summ_{\chi\bmod q} L(\thalf, f\otimes \chi)\overline{L(\thalf, g\otimes \chi)} =
 \begin{cases}
 \frac{2L(1,f\otimes g)}{\zeta(2)} +O(q^{-\delta}) &\text{ if } f\neq g,\\
  \frac{2L(1,\sym^2 f)}{\zeta(2)}\log q +O(q^{-\delta}) &\text{ if } f= g,
 \end{cases}
\end{align}
where $f$ and $g$ are Hecke eigenforms and $\delta>0$, has been obtained by Kowalski, Michel, and Sawin \cite{kms} for prime $q$, though no serious attempt has been made to optimize the error term. We refer the reader to \cite{bfkmms} for applications of this result. 

\subsection*{Acknowledgment} We thank the anonymous referee for their valuable comments.

\section{Sketch.} We briefly indicate the main idea behind the improved error term in Theorem \ref{main} before proceeding to the proof. After using an approximate functional equation and the orthogonality of characters, the error term is roughly
\begin{align}
\label{rough}\max_{1\le NM\le q^{2+\epsilon}}  \Bigg| \frac{1}{\sqrt{MN}}\sum_{\substack{n\equiv m\bmod q \\ n\sim N, m\sim M \\ n\neq m }} \lambda_f(m)d(n)\Bigg| ,
\end{align}
where $\lambda_f(m)$ are the Hecke eigenvalues associated to $f$ and $d(n)$ is the divisor function. This sum is treated using different techniques depending on the the sizes of $N$ and $M$. 

Suppose that $N\ge M$. When $N/M$ is bounded by a reasonably small power of $q$, this is treated by solving a `shifted convolution problem', entailing the spectral theory of automorphic forms. Now suppose $N/M$ is somewhat large, the most difficult case being $N/M\sim q^{\half}$, or equivalently $N\sim q^{\frac32}$ and $M\sim q^\half$ when $NM\sim q^2$. The congruence in \eqref{rough} is detected using additive characters and then Voronoi summation is applied to the $n$ sum (for this step we need a long $n$ sum, hence the assumption $N\ge M$). Writing $N=N_1N_2$ for $N_1\le N_2$, this leads to a sum of the shape (after opening up the divisor function)
\begin{align}
\label{after-poisson} \sum_{\substack{h_1 \sim q/N_1 , h_2\sim q/N_2,  m\sim M}} \lambda_f(m)S(m, h_1 h_2,q).
\end{align}
An effective way to treat the most difficult case $N/M\sim q^{\half}$ is to use Fouvry, Kowalksi, and Michel's \cite{fkm} bounds for bilinear sums of Kloosterman sums, based on their work on algebraic trace functions. To use this, one groups together two of the three variables $h_1,h_2, m$ to form a new variable, which gives a bilinear sum in this new variable and the remaining variable. The bound is non-trivial when one variable is longer than $q^\epsilon$ and one is longer than $q^{\half+\epsilon}$ (the so called P\'olya-Vinogradov range). However this strategy fails when we have $N_1\sim q^\half$ and $N_2\sim q$, for then $h_1\sim M\sim q^\half, h_2\sim 1$, and there is no grouping that satisfies the requirements. This troublesome range was treated by a new bound for bilinear sums of Kloosterman sums obtained by Blomer et al. \cite[Theorem 5.1, equation (5.3)]{bfkmm1}. Later, Shparlinski gave a better bound using his own result on sums of Kloosterman sums \cite[Theorem 2.1]{shp}. The aforementioned bounds for sums of Kloosterman sums are very general because they can be used for arbitrary coefficients in place of $\lambda_f(m)$. On the other hand, we make gains by using the automorphic nature of these coefficients. We transform the sum of Kloosterman sums using Poisson summation to an additively twisted sum of $\lambda_f(m)$, and then employ a Wilton type bound giving square-root cancelation. This strategy yields a bound which essentially forces $N_1$ to be close to $N_2$, thus giving a superior treatment of the critical range $N/M\sim q^{\half}$ and $N_1\sim q^\half$ and $N_2\sim q$ discussed above. Moreover, this approach is simple in that it completely avoids the consideration of bilinear sums of Kloosterman sums. 

The case $M>N$ is somewhat similar after the `switching trick' of Blomer et al. (Voronoi summation) which reduces to sums like \eqref{rough} with $N'\ge M'$ typically, but now without the restriction $N'M'\le q^{2+\epsilon}$. The lack of this restriction in this case weakens the final bound a bit.

\section{Proof of Theorem \ref{main}}

\subsection{Preliminaries} Throughout, $q$ will denote an odd prime number, and $f$ will denote a holomorphic Hecke cusp form or non-holomorphic Hecke-Maass cusp form for $SL_2(\mathbb{Z})$, with Hecke eigenvalues $\lambda_f(n)$.  We will follow the $\epsilon$ convention: that is, $\epsilon$ will always represent an arbitrarily small positive constant, but not necessarily the same one from one occurrence to another.  All implied constants may depend implicitly on $\epsilon$ and $f$. By a `negligible' quantity, we will mean one which is $O(q^{-A})$ for any $A>0$, where the implied constant depends on $A$.  Throughout, we will use $W$, possibly with a subscript, to denote a smooth function compactly supported on the interval $[\half,2]$, with derivatives satisfying
\begin{align*}
W^{(j)}(x)\ll (q^{\epsilon})^j
\end{align*}
for any $j\ge 0$, where the implied constant depends on $j$. Such a weight function may differ from one occurrence to the next even when we use the same notation.  

We will write $e(x)$ to denote $e^{2\pi i x}$. We will use the following elementary orthogonality relation:
\begin{align}
\label{additive} \frac{1}{q}\sum_{a\bmod q} e\Big(\frac{an}{q}\Big) =
\begin{cases}
1 &\text{ if } q|n,\\
0 &\text{ if } q\nmid n.
\end{cases}
\end{align}

We have the average bound \cite[Lemma 1]{iwa3}
\[
\sum_{n\le x} |\lambda_f(n)|^2 \ll x^{1+\epsilon}.
\]
For individual bounds, we have 
\[
|\lambda_f(n)| \ll n^{\theta_f+\epsilon},
\]
where $\theta_f=\frac{7}{64}$ if $f$ is Maass (see \cite{kimsar}) and $\theta_f=0$ if $f$ is holomorphic. 

For smooth sums of Hecke eigenvalues, we note a simple consequence of Voronoi summation: for $M\ge q^\epsilon$, we have
\begin{align}
\label{vor} \sum_{m\ge 1} \lambda_f(m) W\Big(\frac{m}{M}\Big)\ll q^{-A}
\end{align}
for any $A>0$. See for example \cite[Lemmas 2.3, 2.4]{bfkmm1}, where we can take $\vartheta=0$ since the Selberg eigenvalue conjecture is known in level 1.

For (not necessarily smooth) sums of Hecke eigenvalues twisted by additive characters, we have 
\begin{lemma}\label{wilton}
Let $\alpha\in \mathbb{R}$. We have
\[
 \sum_{n\le N} \lambda_f(n) e(n\alpha) \ll N^{\half+\epsilon}.
\]
The implied constant is uniform in $\alpha$, but depends on $f$.
\end{lemma}
\proof
See \cite[Theorem 5.3]{iwa1}, \cite[Theorem 8.1]{iwa2}. 
\endproof

\subsection{Some known bounds}
For $N,M\ge 1$, define
\begin{align}
\label{edef} &E^{\pm}(M,N):= S_1 +S_2 \\
\nonumber&=\frac{1}{\sqrt{MN}} \sum_{\substack{m\equiv \pm n \bmod q\\ m\neq n}} \lambda_f(m) d(n) W_1\Big(\frac{m}{M}\Big) W_2\Big(\frac{n}{N}\Big)
- \frac{1}{q \sqrt{MN} } \sum_{m,n\ge 1}  \lambda_f(m) d(n) W_1\Big(\frac{m}{M}\Big) W_2\Big(\frac{n}{N}\Big).
%\\&=\frac{1}{q\sqrt{MN} } \sum_{\substack{a\bmod q\\ a\neq 0}} \sum_{m,n\ge 1}  \lambda_f(n) d(n) e\Big(\frac{a(m\mp n)}{q}\Big) W_1\Big(\frac{m}{M}\Big) W_2\Big(\frac{n}{N}\Big),
\end{align}
Note that $E^{\pm}(M,N)$ is the same as the quantity $E_{f,E}^{\pm}(M,N) = B_{f,E}^\pm(M,N)$ defined on \cite[page 771]{bfkmm1}. To establish Theorem \ref{main}, we need to prove
\begin{align}
\label{goal} \max_{MN \le q^{2+\epsilon}} |E^{\pm}(M,N)| \ll q^{-\eta_f+\epsilon}.
\end{align}
We begin by recalling some established bounds for $E^{\pm}(M,N)$, but for later use it will be important to not restrict to $MN\le q^{2+\epsilon}$.
 
\begin{lemma}\label{known}
Let $MN\le q^{4+\epsilon}$. We have
\begin{align}
\label{bound1} &E^{\pm}(M,N) \ll 
\begin{cases}
 q^\epsilon \frac{\sqrt{MN}}{q}  &\text{ if } N\ge M\\
 q^\epsilon M^{\theta_f} \frac{\sqrt{MN}}{q} &\text{ if } M> N,
 \end{cases}\\
\label{bound2} &E^{\pm}(M,N) \ll q^\epsilon  \Big(1+\Big(\frac{NM}{q^2}\Big)^\frac14 \Big) \Big(\frac{\max(N,M)}{q\min(N,M)}\Big)^{\frac14} \Big( 1+\Big(\frac{\max(N,M)}{q\min(N,M)}\Big)^{\frac14}\Big).
\end{align}
If $\frac{\max(N,M)}{\min(N,M)}>10$, then
\begin{align}
\label{bound3} E^{\pm}(M,N) \ll q^\epsilon \Big(\frac{q\min(N,M)}{\max(N,M)}\Big)^{\frac12} + q^\epsilon M^{\theta_f}\frac{(NM)^\half }{q^2}.
\end{align}
\end{lemma}
\proof
The first bound \eqref{bound1} is the trivial bound given in \cite[Proposition 3.1]{bfkmm1}.

 The second bound is obtained in \cite[Theorem 3.2]{bfkmm1} by solving a shifted convolution problem using the spectral theory of automorphic forms, with an additional idea required to avoid dependence on the Ramanujan conjecture. However, just reading the statement of that theorem is not enough because it assumes $NM\le q^{2+\epsilon}$. We look at the proof of \cite[Theorem 3.2]{bfkmm1} and find that the bound \eqref{bound2} is contained in the bounds given by \cite[equations (3.12), (3.16)]{bfkmm1}, after dividing by the $(MN)^\half$ factor in \cite[equation (3.4)]{bfkmm1}. This then gives exactly what is seen in the statement of \cite[Theorem 3.2]{bfkmm1}, apart from an extra factor $O(1+(\frac{NM}{q^2})^\frac14)$ which we retain in \eqref{bound2} and which would have been redundant had we been assuming $MN\le q^{2+\epsilon}$, and a missing $O(q^{-\half+\theta_f+\epsilon})$ term, which is specific to the case of $f$ and $g$ both cuspidal (see the bottom of \cite[page 726]{bfkmm1}) and therefore does not appear in \eqref{bound2}.
 
The third bound \eqref{bound3} is obtained in \cite[section 6.4]{bfkmm1} by applying Voronoi summation to both $n$ and $m$ sums and then bounding trivially using Weil's bound for the Kloosterman sum. The assumption $\frac{\max(N,M)}{\min(N,M)}>10$ implies that the condition $n\neq m$ is vacuous in the sum $S_1$ in the definition \eqref{edef} of $E^{\pm}(M,N)$. Then $E^{\pm}(M,N)$ may be expressed as \cite[equation (6.21)]{bfkmm1}, for which we have the bound given in \cite[equation (6.23)]{bfkmm1} plus the trivial bound for the expression in \cite[equation (6.22)]{bfkmm1}. The latter bound gives the term $q^\epsilon M^{\theta_f}\frac{(NM)^\half }{q^2}$ in \eqref{bound3} which does not appear in \cite[equation (6.23)]{bfkmm1} because it is assumed there that $MN\le q^{2+\epsilon}$, while we do not make this assumption. 
\endproof

%%%

\subsection{A new bound} 

Define the sum (closely following the notation in \cite[equation (6.26)]{bfkmm1}) for $M, N_1, N_2\ge 1$ with $N=N_1N_2$ and $N_1\le N_2$:
\[
C^\pm(M,N,N_1,N_2):= \frac{1}{\sqrt{NM}} \sum_{\substack{n_1,n_2,m\ge 1\\ n_1n_2\equiv \pm m \bmod q}}  \lambda_f(m) W_1\Big(\frac{n_1}{N_1}\Big) W_2\Big(\frac{n_2}{N_2}\Big)  W_3\Big(\frac{m}{M}\Big).
\]
We will prove the following simple bound which was not observed in \cite{bfkmm1} or \cite{shp}.
\begin{lemma}\label{new-bound}
Let $ MN< \frac{1}{8} q^3$ and $M\ge q^\epsilon$. We have
\begin{align}
\label{bound5} C^\pm(M,N,N_1,N_2)\ll q^\epsilon \Big(\frac{N_1}{N_2}\Big)^\half + q^{-A},
\end{align}
for any $A>0$.
\end{lemma}
\proof 
Detecting the congruence $n_1n_2\equiv \pm m \bmod q$ using additive characters, we have
\begin{align*}
C^\pm(M,N,N_1,N_2)&= \frac{1}{q\sqrt{NM}} \sum_{a\bmod q} \ \sum_{\substack{n_1,n_2,m\ge 1}}  \lambda_f(m) W_1\Big(\frac{n_1}{N_1}\Big) W_2\Big(\frac{n_2}{N_2}\Big)  W_3\Big(\frac{m}{M}\Big)e\Big(\frac{a(n_1n_2\mp m)}{q}\Big)\\
&=\frac{1}{q\sqrt{NM}} \sum_{\substack{a\bmod q\\ b\bmod q}} \ \sum_{\substack{n_1,n_2,m\ge 1\\ n_2\equiv b\bmod q}}  \lambda_f(m) W_1\Big(\frac{n_1}{N_1}\Big) W_2\Big(\frac{n_2}{N_2}\Big)  W_3\Big(\frac{m}{M}\Big)e\Big(\frac{a(n_1 b \mp m)}{q}\Big).
\end{align*}
Now we apply Poisson summation to the $n_2$ sum (more precisely, we write $n_2=\ell q+b$ and apply Poisson summation to the $\ell$ sum) to get
\begin{multline}
\label{poisson} C^\pm(M,N,N_1,N_2) \\
=\frac{1}{q\sqrt{NM}} \frac{N_2}{q} \sum_{\substack{a\bmod q\\ b\bmod q}} \ \sum_{\substack{n_1,m\ge 1\\ -\infty < k<\infty}}  \lambda_f(m)   e\Big(\frac{bk}{q}\Big)e\Big(\frac{a(n_1 b \mp m)}{q}\Big) W_1\Big(\frac{n_1}{N_1}\Big) \widehat{W_2}\Big(\frac{kN_2}{q}\Big)  W_3\Big(\frac{m}{M}\Big),
\end{multline}
where $\widehat{W_2}(y)=\int_{\mathbb{R}} W_2(x) e(-yx) dx$ is the Fourier transform of $W_2$. By repeated integration by parts, we may restrict the $k$ sum to $|k|< \frac{q^{1+\epsilon}}{N_2}$, up to a negligible error. 

Consider first the contribution to \eqref{poisson} of the terms with $q|n_1$, a condition which implies that $N_1\ge \half q$. This contribution is
\begin{align}
\label{degenerate} \frac{1}{q\sqrt{NM}} \frac{N_2}{q} \sum_{\substack{a\bmod q\\ b\bmod q}} \ \sum_{\substack{n_1,m\ge 1 \\ 0\le |k| < \frac{q^{1+\epsilon}}{N_2} \\ q|n_1}}  \lambda_f(m) e\Big(\frac{bk}{q}\Big) e\Big(\frac{a( \mp m)}{q}\Big) W_1\Big(\frac{n_1}{N_1}\Big) \widehat{W_2}\Big(\frac{kN_2}{q}\Big)  W_3\Big(\frac{m}{M}\Big) + O(q^{-A})
\end{align}
for any $A>0$. Using \eqref{additive}, we see that the $a$ sum vanishes unless $q$ divides $m$, in which case we must have $M\ge \half q$. But since we are assuming $MN< \frac{1}{8} q^3$, we cannot have $M\ge \half q$ and $N_2\ge  N_1\ge \half q$, so it follows that \eqref{degenerate} is $O(q^{-A})$. 

It remains to consider the contribution to \eqref{poisson} of the terms with $q\nmid n_1$. Evaluating the $b$ sum using \eqref{additive} enforces $a \equiv k \overline{n_1} \bmod q$, where $\overline{n_1}$ denotes the multiplicative inverse of $n_1$ modulo $q$. Thus we have that \eqref{poisson} equals
\[
\frac{N_2}{q\sqrt{NM}}  \sum_{\substack{n_1,m\ge 1 \\ 0\le  |k| < \frac{q^{1+\epsilon}}{N_2} \\ q\nmid n_1}}   \lambda_f(m)   e\Big(\frac{ \mp  m k\overline{n_1} )}{q}\Big) W_1\Big(\frac{n_1}{N_1}\Big) \widehat{W_2}\Big(\frac{kN_2}{q}\Big)  W_3\Big(\frac{m}{M}\Big) + O(q^{-A}).
\]
The contribution to the sum of the term $k=0$ is negligible by using \eqref{vor}. We estimate the contribution of the terms $k\neq 0$ by bounding the $m$ sum using Lemma \ref{wilton} and the rest trivially, to get
\[
\frac{N_2}{q\sqrt{NM}}  \frac{q^{1+\epsilon}}{N_2}  N_1 \sqrt{M} \ll q^\epsilon \Big(\frac{N_1}{N_2}\Big)^\half.
\]
\endproof
In order to use the new bound towards our goal \eqref{goal}, we need
\begin{lemma} \label{lemma-use}
If $\frac{\max(N,M)}{\min(N,M)}>10$ and $M\ge q^\epsilon$, then
\begin{align}
\label{use0} |E^{\pm}(M,N)| \ll \max_{\substack{N=N_1N_2\\ N_1\le N_2}} |C^\pm(M,N,N_1,N_2)| + O(q^{-A}).
\end{align}
for any $A>0$.
\end{lemma}
\proof
The assumption $M\ge q^\epsilon$ implies by \eqref{vor} that the sum $S_2$ in the definition \eqref{edef} of $E^{\pm}(M,N)$ is negligible. The assumption $\frac{\max(N,M)}{\min(N,M)}>10$ implies that the condition $n\neq m$ is vacuous in the sum $S_1$ in the definition \eqref{edef} of $E^{\pm}(M,N)$.  After removing this condition in $S_1$, we may open the divisor function and apply smooth partitions of unity, as explained in \cite[section 6.4.1]{bfkmm1}, to obtain \eqref{use0}.
\endproof

%%%

\subsection{The case $M>N$} \label{m>n case}

Recall our goal \eqref{goal}. We need to start with the `switching trick' of \cite{bfkmm1}.

\begin{lemma}\label{use}
Let $MN\le q^{2+\epsilon}$ and $M>N$. Define
\[
M^*=\frac{q^2}{M}, \ \ N^*=\frac{q^2}{N}.
\]
We have
\begin{align}
\label{switch2} E^\pm(M,N) \ll  q^\epsilon \max_{\substack{ 1\le N'\le q^\epsilon N^* \\ 1\le M'\le  q^\epsilon M^* }} | E^\pm(M',N') | + O(q^{-\frac16+\epsilon}).
\end{align}
\end{lemma}
\proof
By applying Voronoi summation to the $n$ and $m$ sums in \eqref{edef}, Blomer et al. \cite[page 765]{bfkmm1} showed\footnote{They actually have an additional factor $(M^*/M')^{2\theta_f}$, where $\theta_f=7/64$, which arises from the use of their Lemma 2.4, but this factor is not needed because we work with Hecke eigenforms of level 1. In this setting, the Selberg eigenvalue conjecture is known, so that we may take $\vartheta=0$ in their Lemma 2.4.}
\begin{align}
\label{switch} E^\pm(M,N) &\ll q^\epsilon \max_{\substack{ 1\le N'\le q^\epsilon N^* \\ 1\le M'\le  q^\epsilon M^* }} \Bigg| \frac{1}{\sqrt{M^* N^* }}  \sum_{\substack{m,n\ge1\\ m\equiv\pm n \bmod q}} \lambda_f(m)d(n) W_1\Big(\frac{m}{M'}\Big)W_2\Big(\frac{n}{N'}\Big) \Bigg| \\
\nonumber &+  q^\epsilon \max_{\substack{ 1\le N'\le q^\epsilon N^* \\ 1\le M'\le  q^\epsilon M^* }} \Bigg| \frac{1 }{q \sqrt{M^*N^*} } \sum_{m,n\ge 1}  \lambda_f(m) d(n) W_1\Big(\frac{m}{M'}\Big) W_2\Big(\frac{n}{N'}\Big)\Bigg| + O(q^{-1+\epsilon}),
\end{align}
where
\[
M^*=\frac{q^2}{M}, \ \ N^*=\frac{q^2}{N}.
\]
The contribution of the diagonal terms $n=m$ in the first sum on the right hand side of \eqref{switch} is bounded absolutely by 
\[ 
q^\epsilon \frac{1 }{\sqrt{M^* N^* }}  M' \ll  q^\epsilon \frac{1 }{\sqrt{M^* N^* }}  M^*   \ll  q^\epsilon \sqrt\frac{M^*}{N^*} = q^\epsilon \sqrt\frac{N}{M}.
\]
We can assume this is $O(q^{-\frac16+\epsilon})$ or else by \eqref{bound2} we have $E^\pm(M,N)\ll q^{-\frac16+\epsilon}$. Thus at the cost of an admissible error term, the diagonal terms can be removed. By \eqref{vor}, the second sum on the right hand side of \eqref{switch} is negligible unless $M'<q^\epsilon$. Thus this sum is trivially bounded by
\[
q^\epsilon \frac{1}{q\sqrt{M^* N^* }} N^* \ll \frac{q^\epsilon}{q} \sqrt\frac{M}{N}.
\]
We can assume this is $O(q^{-\frac14+\epsilon})$ or else by \eqref{bound3} we have $E^\pm(M,N)\ll q^{-\frac14+\epsilon}$. The lemma follows.
\endproof

We will apply Lemma \ref{new-bound} and Lemma \ref{lemma-use} to bound $E^\pm(M',N')$. For this, we need to verify the conditions of these lemmas. Suppose that the condition $\frac{\max(N',M')}{\min(N',M')}>10$ does not hold. Then by \eqref{bound1} and \eqref{bound2}, we would have
\begin{multline*}
E^\pm(M',N')
\ll q^\epsilon\min\Bigg( (M')^{\theta_f} \Big(\frac{N'M'}{q^2}\Big)^{-\half}, \Big(\frac{N'M'}{q^2}\Big)^{\frac14}q^{-\frac14}  \Bigg)\\
 \ll q^\epsilon\min\Bigg( q^{\frac{7}{32}} \Big(\frac{N'M'}{q^2}\Big)^{-\half},  \Big(\frac{N'M'}{q^2}\Big)^{\frac14} q^{-\frac14}  \Bigg)\ll q^{-\frac{3}{32}+\epsilon},
\end{multline*}
and we would be done by \eqref{switch2}. So we can assume that $\frac{\max(N',M')}{\min(N',M')}>10$ holds. Now suppose that the condition $M'\ge q^\epsilon$ does not hold. Then by \eqref{bound1} and \eqref{bound3}, we would have
\[
E^\pm(M',N')\ll q^\epsilon\min\Bigg(  \frac{\sqrt{N'}}{q}, \Big(\frac{q}{N'}\Big)^\half + \frac{\sqrt{N'}}{q^2} \Bigg) \ll  \Bigg(  \frac{\sqrt{N'}}{q}, \Big(\frac{q}{N'}\Big)^\half + q^{-1}\Bigg) \ll  q^{-\frac14+\epsilon},
\]
and we would be done by \eqref{switch2}. So we can assume that $M'\ge q^\epsilon$ holds. Finally suppose that the condition $M'N'<\frac18 q^3$ does not hold. Then by \eqref{bound1}, we would have $M^*N^*\gg q^3$, which implies $MN\ll q$. Then by \eqref{bound1}, we would have $E^\pm(M',N')\ll (M')^{\theta_f} q^{-\half+\epsilon}\ll q^{\frac{7}{32}-\half+\epsilon}$ and we would be done by \eqref{switch2}. So we can assume the condition $M'N'<\frac18 q^3$ as well.

Now we are ready to prove
\begin{lemma}\label{m>n} We have
\[
 \max_{\substack{MN \le q^{2+\epsilon}\\ M> N}} |E^{\pm}(M,N)| \ll q^{-\eta_f+\epsilon}.
 \]
\end{lemma}
\proof
We use Lemma \ref{new-bound} and Lemma \ref{lemma-use} to bound $E^{\pm}(M',N')$, which implies a bound for $E^{\pm}(M,N)$ via Lemma \ref{use}. We add our bound to the collection of bounds given in \cite{bfkmm1}. We need to find the worst case (maximum) of the minimum of all these bounds for $MN \le q^{2+\epsilon}$ and $M>N$. This leads to a linear optimization problem to maximize the minimum of the exponents of $q$ from each bound. Such a problem was solved in \cite{bfkmm1} by a computer search, and we do the same by simply inserting the exponent of $q$ from our bound for $C^\pm(M',N',N_1,N_2)$, given in \eqref{bound5}, into the existing Mathematica code from \cite[section 7.4.2]{bfkmm1}. The underlined portion below indicates the new addition to this code. This returns an exponent of $-\frac{5}{152}$ for $f$ Maass and $-\frac{1}{22}$ for $f$ holomorphic (in the latter case we also omit the term $7m/64$ from the first entry in the existing code, corresponding to the bound \eqref{bound1}, since $\theta_f=0$).

\

\

\noindent \begin{tabular}{ll}
{\tt In[1] := } &   {\tt Maximize[\{Min[(m + n)/2 - 1 + 7m/64, 
   Max[(m - n - 1)/2, (m - n - 1)/4], }\\ & {\tt (1 + n - m)/2, 
   7/64(mstar-mprime) + Max[(2 mprime - mstar - nstar + 1)/
     2, }\\ & {\tt (2 n1 + 2 mprime - 2 - mstar - nstar)/2], \underline{(n1 - n2)/2,}
   }\\ & {\tt 7/64(mstar-mprime) + Max[(mprime + n2 + 1/2 - nstar - mstar)/
     2, }\\ & {\tt (2 mprime + n1 - mstar - nstar)/
     2, (2 n1 + 2 mprime - 2 - mstar - nstar)/2], }\\ & {\tt 
   7/64(mstar-mprime) + mprime + n1 - (mstar + nstar)/2, }\\
  & {\tt 
   7/64(mstar-mprime) + Max[Min[2 n1 - (mstar + nstar)/2, }\\ & {\tt 
     2/3 nprime + n1 + 1/2 - 1/6 mprime - 
      n2 - (mstar + nstar)/2], }\\ & {\tt (2 mprime - mstar - nstar)/
     2, (2 mprime + 2 n1 - 2 - mstar - nstar)/2, }\\ & {\tt 
    2 mprime - n2 - (mstar + nstar)/2]], 0 <= n, n <= m, m + n <= 2, 
}\\ & {\tt   nstar == 2 - n, mstar == 2 - m, 0 <= nprime, nprime  <= nstar, 
 }\\ & {\tt  0 <= mprime, mprime <= mstar, n1 + n2 == nprime, 0 <= n1, 
  n1 <= n2\},}\\ & {\tt  \{m, n, n1, n2, nprime, nstar, mprime, mstar\}]}\\[0.3cm]
 {\tt Out[1] := } & $\Bigl\{ -\frac{5}{152}, \Bigl\{ {\tt m} \rightarrow \frac{24}{19}, \,\,{\tt n} \rightarrow \frac{15}{38}, \,\, {\tt n1} \rightarrow \frac{61}{76}, \,\, {\tt n2} \rightarrow \frac{61}{76}, \,\, {\tt nprime} \rightarrow \frac{61}{38}, \,\, {\tt nstar} \rightarrow \frac{61}{38}$,\\ & ${\tt mprime} \rightarrow \frac{14}{19}, \,\, {\tt mstar} \rightarrow \frac{14}{19}
  \Bigr\} \Bigr\}$
  \end{tabular} 

\

\

\noindent \begin{tabular}{ll}
{\tt In[2] := } &   {\tt Maximize[\{Min[(m + n)/2 - 1, 
   Max[(m - n - 1)/2, (m - n - 1)/4], }\\ & {\tt (1 + n - m)/2, 
   7/64(mstar-mprime) + Max[(2 mprime - mstar - nstar + 1)/
     2, }\\ & {\tt (2 n1 + 2 mprime - 2 - mstar - nstar)/2], \underline{(n1 - n2)/2,}
   }\\ & {\tt 7/64(mstar-mprime) + Max[(mprime + n2 + 1/2 - nstar - mstar)/
     2, }\\ & {\tt (2 mprime + n1 - mstar - nstar)/
     2, (2 n1 + 2 mprime - 2 - mstar - nstar)/2], }\\ & {\tt 
   7/64(mstar-mprime) + mprime + n1 - (mstar + nstar)/2, }\\
  & {\tt 
   7/64(mstar-mprime) + Max[Min[2 n1 - (mstar + nstar)/2, }\\ & {\tt 
     2/3 nprime + n1 + 1/2 - 1/6 mprime - 
      n2 - (mstar + nstar)/2], }\\ & {\tt (2 mprime - mstar - nstar)/
     2, (2 mprime + 2 n1 - 2 - mstar - nstar)/2, }\\ & {\tt 
    2 mprime - n2 - (mstar + nstar)/2]], 0 <= n, n <= m, m + n <= 2, 
}\\ & {\tt   nstar == 2 - n, mstar == 2 - m, 0 <= nprime, nprime  <= nstar, 
 }\\ & {\tt  0 <= mprime, mprime <= mstar, n1 + n2 == nprime, 0 <= n1, 
  n1 <= n2\},}\\ & {\tt  \{m, n, n1, n2, nprime, nstar, mprime, mstar\}]}\\[0.3cm]
 {\tt Out[2] := } & $\Bigl\{ -\frac{1}{22}, \Bigl\{ {\tt m} \rightarrow \frac{15}{11}, \,\,{\tt n} \rightarrow \frac{6}{11}, \,\, {\tt n1} \rightarrow \frac{8}{11}, \,\, {\tt n2} \rightarrow \frac{8}{11}, \,\, {\tt nprime} \rightarrow \frac{16}{11}, \,\, {\tt nstar} \rightarrow \frac{16}{11}$,\\ & ${\tt mprime} \rightarrow \frac{7}{11}, \,\, {\tt mstar} \rightarrow \frac{7}{11}
  \Bigr\} \Bigr\}$
  \end{tabular} 
  
  \endproof
  
%%%

\subsection{The case $N\ge M$} Recall our goal \eqref{goal}. Assuming $MN\le q^{2+\epsilon}$ and $N\ge M$, we will apply Lemma \ref{new-bound} and Lemma \ref{lemma-use} to bound $E^\pm(M,N)$. For this purpose, we need to verify the conditions $\frac{\max(N,M)}{\min(N,M)}>10, M\ge q^\epsilon,$ and $MN<\frac18 q^3$ of these lemmas. This is done just as in section \ref{m>n case}  and is in fact easier. This is because the factor $M^{\theta_f}$ does not occur when using \eqref{bound1} because we assume that $N\ge M$, and the condition $MN<\frac18 q^3$ is automatic for $q$ large enough by the assumption $MN\le q^{2+\epsilon}$.

Now we prove the following, which together with Lemma \ref{m>n} gives Theorem \ref{main}.
\begin{lemma} We have
\[
 \max_{\substack{MN \le q^{2+\epsilon}\\ N\ge M}} |E^{\pm}(M,N)| \ll q^{-\frac{1}{20}+\epsilon}.
 \]
\end{lemma}
\proof
As before, we simply insert our bound from Lemma \ref{new-bound} into the existing Mathematica code, this time from \cite[section 7.4.1]{bfkmm1}, and compute.

\

\

\noindent \begin{tabular}{ll}
{\tt In[3] := } &   {\tt Maximize[\{Min[(m + n)/2 - 1, 
   Max[(n - m - 1)/2, (n - m - 1)/4],}\\ & {\tt (1 + m - n)/2, 
   Max[(m + 1 - n)/2, (2 n1 + m - 2 - n)/2], }\\ 
   & {\tt Max[1/4 - n1/2, (m - n2)/2, (2 n1 + m - n - 2)/2], (m + 2 n1 - n)/
    2, }\\ & {\tt Max[Min[2 n1 - (m + n)/2, 
     n/6 + n1 + 1/2 - n2 - 2 m/3], (m - n)/2, }\\ & {\tt (m + 2 n1 - 2 - n)/2, 
    3 m/2 - n2 - n/2], 
   If[m + n1prime <= 
     1, }\\ & {\tt (2 m + 2 n1prime + 2 n2prime - 1 - n1circ - n2circ - m)/2 
   }\\ & {\tt +  Max[-n2prime/2, 1/4 - (m + n1prime)/2], 10], \underline{(n1 - n2)/2,} 
 }\\ & {\tt  If[m + n2prime <= 
     2 n1prime, (2 m + 2 n1prime + 2 n2prime - 1 }\\ & {\tt - n1circ - n2circ - 
        m)/2 + 1/4 - 5 n1prime/12 - (m + n2prime)/6, 10]], }\\ & {\tt m >= 0, 
  n >= m, m + n <= 2, n1 + n2 == n, n1 <= n2, n1 >= 0, }\\ & {\tt  n1prime >= 0, 
  n1prime <= n1circ, n1circ == 1 - n1,  n2prime >= 0, }\\ & {\tt
  n2prime <= n2circ, n2circ == 1 - n2\}, \{m, n, n1, n2, n1prime, }\\ & {\tt
  n2prime, n1circ, n2circ\}]}\\[0.3cm]
 {\tt Out[3] := } & $\Bigl\{ -\frac{1}{20}, \Bigl\{ {\tt m} \rightarrow \frac{3}{5}, \,\,{\tt n} \rightarrow \frac{7}{5}, \,\, {\tt n1} \rightarrow \frac{7}{10}, \,\, {\tt n2} \rightarrow \frac{7}{10}, \,\, {\tt n1prime} \rightarrow \frac{3}{10}, \,\, {\tt n1prime} \rightarrow \frac{3}{10}$,\\ & ${\tt n1circ} \rightarrow \frac{3}{10}, \,\, {\tt n1circ} \rightarrow \frac{3}{10}
  \Bigr\} \Bigr\}$
  \end{tabular} 
  
  \endproof

%%%

\bibliographystyle{amsplain} 

\bibliography{mixed-moment}

\end{document}